\documentclass[10pt]{article}
\usepackage{cite}
\usepackage{mathrsfs}
\usepackage{amsfonts}
\usepackage{amsmath}
\usepackage{amsfonts,amssymb}
\usepackage{dsfont}
\usepackage{curves}
\usepackage{mathrsfs}
\usepackage{pifont}
\usepackage{amssymb}
\usepackage{graphicx}
\allowdisplaybreaks

\numberwithin{equation}{section}

\date{}

\begin{document}
\title{Some existence theorems on path-factor critical avoidable graphs
%\thanks{Supported by the National Natural Science Foundation of China (Grant No. 11371009, 11501256,
%61503160), the National Social Science Foundation of China (Grant No. 14AGL001), and sponsored by 333
%Project of Jiangsu Province.}
}
\author{\small  Sizhong Zhou$^{1}$\footnote{Corresponding
author. E-mail address: zsz\_cumt@163.com (S. Zhou)}, Hongxia Liu$^{2}$\\
\small  1. School of Science, Jiangsu University of Science and Technology,\\
\small  Zhenjiang, Jiangsu 212100, China\\
\small  2. School of Mathematics and Information Sciences, Yantai University,\\
\small  Yantai, Shandong 264005, China\\
}

\maketitle
\begin{abstract}
\noindent A spanning subgraph $F$ of $G$ is called a path factor if every component of $F$ is a path of order at least 2. Let $k\geq2$ be an
integer. A $P_{\geq k}$-factor of $G$ means a path factor in which every component has at least $k$ vertices. A graph $G$ is called a
$P_{\geq k}$-factor avoidable graph if for any $e\in E(G)$, $G$ has a $P_{\geq k}$-factor avoiding $e$. A graph $G$ is called a
$(P_{\geq k},n)$-factor critical avoidable graph if for any $W\subseteq V(G)$ with $|W|=n$, $G-W$ is a $P_{\geq k}$-factor avoidable graph.
In other words, $G$ is $(P_{\geq k},n)$-factor critical avoidable if for any $W\subseteq V(G)$ with $|W|=n$ and any $e\in E(G-W)$, $G-W-e$
admits a $P_{\geq k}$-factor. In this article, we  verify that (\romannumeral1) an $(n+r+2)$-connected graph $G$ is $(P_{\geq2},n)$-factor
critical avoidable if $I(G)>\frac{n+r+3}{2(r+2)}$; (\romannumeral2) an $(n+r+2)$-connected graph $G$ is $(P_{\geq3},n)$-factor critical
avoidable if $t(G)>\frac{n+r+2}{2(r+2)}$; (\romannumeral3) an $(n+r+2)$-connected graph $G$ is $(P_{\geq3},n)$-factor critical avoidable
if $I(G)>\frac{n+3(r+2)}{2(r+2)}$; where $n$ and $r$ are two nonnegative integers.
\\
\begin{flushleft}

(2020) Mathematics Subject Classification: 05C70, 05C38, 90B10
\end{flushleft}
\end{abstract}

\section{Introduction}

In this work, we discuss only finite, undirected and simple graphs. We denote by $G=(V(G),E(G))$ a graph, where $V(G)$ denotes the vertex
set of $G$ and $E(G)$ denotes the edge set of $G$. For a vertex $x$ of $G$, the degree of $x$ in $G$, denoted by $d_G(x)$, is the number of
vertices adjacent to $x$ in $G$. For a vertex subset $X$ of $G$, $G[X]$ denotes the subgraph of $G$ induced by $X$, and $G-X$ denotes the
subgraph derived from $G$ by removing all vertices in $X$. For an edge subset $E'$ of $G$, $G-E'$ denotes the subgraph acquired from $G$
by deleting all edges in $E'$. For a vertex (or an edge) subset $Q$, we denote $G-Q$ by $G-u$ for convenience if $Q=\{u\}$. Let $i(G)$,
$\omega(G)$ and $\kappa(G)$ denote the number of isolated vertices, the number of connected components and the vertex connectivity of $G$,
respectively. We use $K_n$ and $P_n$ to denote the complete graph and the path with $n$ vertices, respectively. Let $G_1$ and $G_2$ be two
graphs. Then the join $G_1+G_2$ denotes the graph with vertex set $V(G_1+G_2)=V(G_1)\cup V(G_2)$ and edge set
$$
E(G_1+G_2)=E(G_1)\cup E(G_2)\cup\{uv:u\in V(G_1),v\in V(G_2)\}.
$$

The toughness of a graph $G$, denoted by $t(G)$, was first introduced by Chv\'atal \cite{C}. If $G$ is not complete, then
$$
t(G)=\min\left\{\frac{|X|}{\omega(G-X)}:X\subseteq V(G),\omega(G-X)\geq2\right\};
$$
otherwise, $t(G)=+\infty$.

The isolated toughness of a graph $G$, denoted by $I(G)$, was first introduced by Yang, Ma and Liu \cite{YML}. If $G$ is not complete, then
$$
I(G)=\min\left\{\frac{|X|}{i(G-X)}:X\subseteq V(G),i(G-X)\geq2\right\};
$$
otherwise, $I(G)=+\infty$.

A spanning subgraph $F$ of $G$ is called a path factor if every component of $F$ is a path of order at least 2. Let $k\geq2$ be an integer. A
$P_{\geq k}$-factor of $G$ means a path factor in which every component has at least $k$ vertices.

Las Vergnas \cite{V} showed a necessary and sufficient condition for graphs to possess $P_{\geq2}$-factors.

\medskip

\noindent{\textbf{Theorem 1}} (\cite{V}). A graph $G$ possesses a $P_{\geq2}$-factor if and only if $G$ satisfies
$$
i(G-X)\leq2|X|
$$
for every vertex subset $X$ of $G$.

\medskip

A graph $H$ is factor-critical if any induced subgraph with $|V(H)|-1$ vertices has a perfect matching. A graph $R$ is called a sun if $R=K_1$,
$R=K_2$ or $R$ is the corona of a factor-critical graph $H$ with at least three vertices, namely, $R$ is acquired from $H$ by adding a new
vertex $z=z(y)$ together with a new edge $yz$ for any $y\in V(H)$ to $H$ (Figure 1, which was shown by Kano, Lu and Yu \cite{KLY}). We easily
see that $d_R(z)=1$. In particular, a sun with at least six vertices is called a big sun. Let $sun(G)$ denote the number of sun components of
$G$. In fact, $i(G)\leq sun(G)\leq\omega(G)$.

\begin{figure}
  \centering
  \includegraphics[height=1.6in,width=3in]{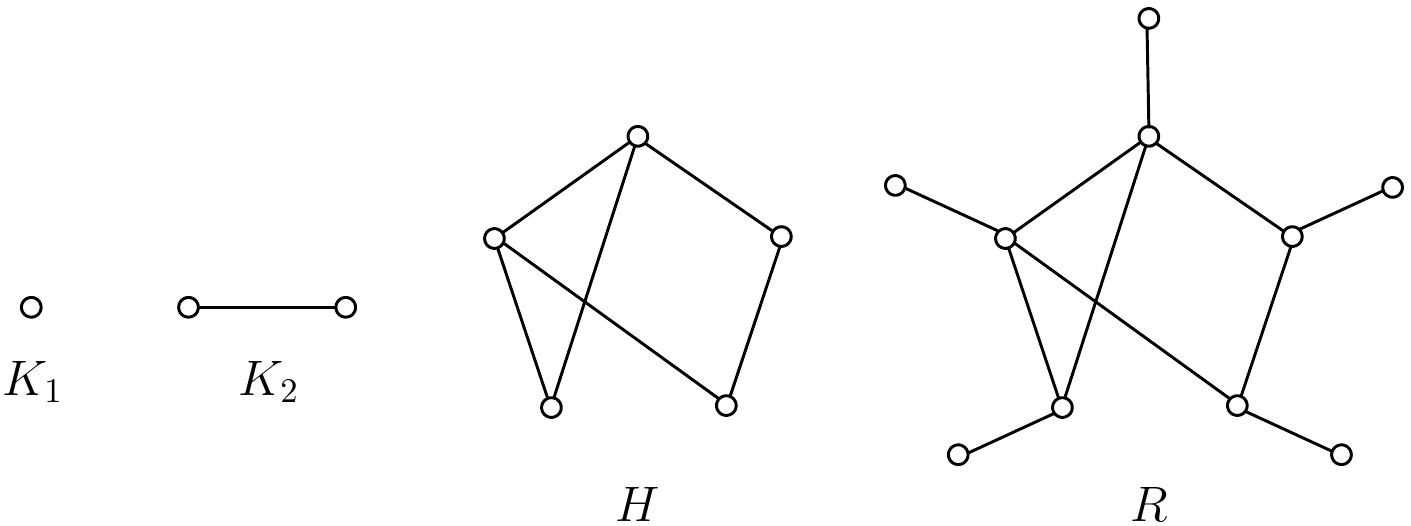}\\
  \caption{A factor-critical graph $H$ and the sun $R$ obtained from $H$.}\label{11}
\end{figure}

Kaneko \cite{K} gave a necessary and sufficient condition for graphs admitting $P_{\geq3}$-factors. Kano, Katona and Kir\'aly \cite{KKK} gave
a simple proof.

\medskip

\noindent{\textbf{Theorem 2}} (\cite{K,KKK}). A graph $G$ contains a $P_{\geq3}$-factor if and only if $G$ satisfies
$$
sun(G-X)\leq2|X|
$$
for every vertex subset $X$ of $G$.

\medskip

In recent years, many results on path factors were derived. Kelmans \cite{Kp} raised some results on the existence of path factors in claw-free
graphs. Ando et al \cite{AEKKM} derived a minimum degree condition for a claw-free graph to have a path factor. Kano, Lee and Suzuki \cite{KLS}
verified that every connected cubic bipartite graph with at least eight vertices admits a $P_{\geq8}$-factor. Egawa and Furuya \cite{EF} showed
some sufficient conditions for graphs to have path factors. Kano, Lu and Yu \cite{KLY} presented a sufficient condition for the existence of
$P_{\geq3}$-factor. Wu \cite{Wp}, Zhou et al \cite{Zp,ZBP,ZSL,ZWX,ZWB} derived some sufficient conditions for graphs to possess $P_{\geq3}$-factors with given
properties. Gao, Wang and Chen \cite{GWC} posed some tight bounds for the existence of $P_{\geq3}$-factors in graphs. Dauer, Katona, Kratsch
and Veldman \cite{BKKV}, Gao, Guirao and Chen  \cite{GGC}, Liu and Zhang \cite{LZ} established some relationships between toughness and graph
factors. Gao and Wang \cite{GW}, Gao, Liang and Chen \cite{GLC} established some relationships between isolated toughness and graph factors.
More results on graph factors were acquired by Zhou \cite{Za,Zb,Zha,Zd,Zr}, Zhou and Liu \cite{ZLt}, Zhou, Xu and Sun
\cite{ZXS}, Wang and Zhang \cite{WZo,WZr}, Yuan and Hao \cite{YHa}.

A graph $G$ is called a $P_{\geq k}$-factor avoidable graph if for any $e\in E(G)$, $G$ has a $P_{\geq k}$-factor avoiding $e$. A graph $G$ is
called a $(P_{\geq k},n)$-factor critical avoidable graph if for any $W\subseteq V(G)$ with $|W|=n$, $G-W$ is a $P_{\geq k}$-factor avoidable
graph. In other words, $G$ is $(P_{\geq k},n)$-factor critical avoidable if for any $W\subseteq V(G)$ with $|W|=n$ and any $e\in E(G-W)$,
$G-W-e$ contains a $P_{\geq k}$-factor.

Zhou \cite{Zhs} acquired some toughness or isolated toughness conditions for graphs to be $(P_{\geq k},n)$-factor critical avoidable graphs for
$k=2,3$.

\medskip

\noindent{\textbf{Theorem 3}} (\cite{Zhs}). An $(n+2)$-connected graph $G$ is $(P_{\geq 2},n)$-factor critical avoidable if its isolated toughness
$I(G)>\frac{n+2}{3}$, where $n\geq0$ is an integer.

\medskip

\noindent{\textbf{Theorem 4}} (\cite{Zhs}). An $(n+2)$-connected graph $G$ is $(P_{\geq 3},n)$-factor critical avoidable if its toughness
$t(G)>\frac{n+1}{2}$, where $n\geq0$ is an integer.

\medskip

\noindent{\textbf{Theorem 5}} (\cite{Zhs}). An $(n+2)$-connected graph $G$ is $(P_{\geq 3},n)$-factor critical avoidable if its isolated toughness
$I(G)>\frac{n+3}{2}$, where $n\geq0$ is an integer.

\medskip

In this article, we proceed to investigate $(P_{\geq k},n)$-factor critical avoidable graphs and derive three more general results on
$(P_{\geq k},n)$-factor critical avoidable graphs depending on toughness and isolated toughness, which are shown in Sections 2 and 3.

\section{$(P_{\geq2},n)$-factor critical avoidable graphs}

In this section, we pose a sufficient conditions using isolated toughness for graphs to be $(P_{\geq2},n)$-factor critical
avoidable graphs, which is an improvements of Theorem 3 for $n\geq1$.

\medskip

\noindent{\textbf{Theorem 6.}} Let $n$ and $r$ be two nonnegative integers, and let $G$ be an $(n+r+2)$-connected graph. If its isolated
toughness $I(G)>\frac{n+r+3}{2(r+2)}$, then $G$ is $(P_{\geq2},n)$-factor critical avoidable.

\medskip

\noindent{\it Proof.}  Obviously, Theorem 6 holds for a complete graph. Next, we assume that $G$ is not a complete graph. Let $H=G-W-e$
for any $W\subseteq V(G)$ with $|W|=n$ and any $e\in E(G-W)$. It suffices to claim that $H$ has a $P_{\geq2}$-factor. Suppose that $H$
has no $P_{\geq2}$-factor. Then it follows from Theorem 1 that
$$
i(H-X)\geq2|X|+1\eqno(1)
$$
for some $X\subseteq V(H)$.

Since $G$ is $(n+r+2)$-connected, $H$ is $(r+1)$-connected.

\noindent{\bf Claim 1.} $|X|\geq r+2$.

\noindent{\it Proof.} If $|X|=0$, then it follows from (1) that $i(H)\geq1$, which contradicts that $H$ is $(r+1)$-connected. Next, we
discuss $1\leq|X|\leq r+1$.

According to (1), we derive $i(H-X)\geq2|X|+1\geq3$. Thus, we have
$$
i(G-W-X)\geq i(G-W-X-e)-2=i(H-X)-2\geq3-2=1,
$$
which implies that there exists an isolated vertex $v$ in $G-W-X$, namely, $d_{G-W-X}(v)=0$. Combining this with $1\leq|X|\leq r+1$, we
get
$$
d_G(v)\leq d_{G-W-X}(v)+|W|+|X|=0+n+|X|\leq n+r+1,
$$
which contradicts that $G$ is $(n+r+2)$-connected. Hence, $|X|\geq r+2$. Claim 1 is proved. \hfill $\Box$

In fact, $i(G-W-X-e)\geq i(G-W-X)\geq i(G-W-X-e)-2$. The following proof is divided into three cases.

\noindent{\bf Case 1.} $i(G-W-X)=i(G-W-X-e)-2$.

In this case, it is obvious that there exists a $K_2$ component in $G-W-X$ with $e\in E(K_2)$. Let $u\in V(K_2)$. Then by (1) and Claim 1, we
deduce
\begin{eqnarray*}
i(G-W-X-u)&=&i(G-W-X)+1=i(G-W-X-e)-2+1\\
&=&i(H-X)-1\geq2|X|+1-1=2|X|\\
&\geq&2(r+2)\geq4,
\end{eqnarray*}
and so
$$
I(G)\leq\frac{|W\cup X\cup\{u\}|}{i(G-W-X-u)}\leq\frac{n+|X|+1}{2|X|}=\frac{1}{2}+\frac{n+1}{2|X|}\leq\frac{1}{2}+\frac{n+1}{2(r+2)}=\frac{n+r+3}{2(r+2)},
$$
which contradicts $I(G)>\frac{n+r+3}{2(r+2)}$.

\noindent{\bf Case 2.} $i(G-W-X)=i(G-W-X-e)-1$.

In this case, there exists a vertex $u$ such that $d_{G-W-X}(u)=1$. Let $v$ be an unique vertex adjacent to $u$ in $G-W-X$, and $e=uv$. Similar to
this discussion of Case 1, we easily deduce
$$
i(G-W-X-v)\geq2|X|+1=2(r+2)+1\geq5,
$$
and so
\begin{eqnarray*}
I(G)&\leq&\frac{|W\cup X\cup\{v\}|}{i(G-W-X-v)}\leq\frac{n+|X|+1}{2|X|+1}=\frac{1}{2}+\frac{n+\frac{1}{2}}{2|X|+1}\\
&<&\frac{1}{2}+\frac{n+1}{2|X|}\leq\frac{1}{2}+\frac{n+1}{2(r+2)}=\frac{n+r+3}{2(r+2)},
\end{eqnarray*}
which contradicts $I(G)>\frac{n+r+3}{2(r+2)}$.

\noindent{\bf Case 3.} $i(G-W-X)=i(G-W-X-e)$.

According to (1), we obtain
$$
I(G)\leq\frac{|W\cup X|}{i(G-W-X)}=\frac{n+|X|}{i(G-W-X-e)}=\frac{n+|X|}{i(H-X)}\leq\frac{n+|X|}{2|X|+1}.\eqno(2)
$$

If $n=0$, then from (2) we have
$$
I(G)\leq\frac{|X|}{2|X|+1}<\frac{1}{2},
$$
which contradicts $I(G)>\frac{r+3}{2(r+2)}>\frac{1}{2}$.

If $n\geq1$, then it follows from (2) and Claim 1 that
$$
I(G)\leq\frac{n+|X|}{2|X|+1}=\frac{1}{2}+\frac{n-\frac{1}{2}}{2|X|+1}<\frac{1}{2}+\frac{n}{2|X|}\leq\frac{1}{2}+\frac{n}{2(r+2)}=\frac{n+r+2}{2(r+2)},
$$
which contradicts $I(G)>\frac{n+r+3}{2(r+2)}$. This completes the proof of Theorem 6. \hfill $\Box$

\medskip

\noindent{\bf Remark 1.} Next, we show that the condition on $I(G)$ in Theorem 6 is sharp.

Let $G=K_{n+r+2}+((2r+3)K_1\cup K_2)$, where $n$ and $r$ are two nonnegative integers with $n\geq r+1$. Clearly, $G$ is $(n+r+2)$-connected and
$I(G)=\frac{n+r+3}{2(r+2)}$. Let $W\subseteq V(K_{n+r+2})\subseteq V(G)$ with $|W|=n$ and $e\in E(K_2)$. Then $G-W-e=K_{r+2}+((2r+5)K_1)$.
Let $X=V(K_{r+2})\subseteq V(G-W-e)$. Then we possess
$$
i(G-W-e-X)=2r+5>2(r+2)=2|X|.
$$
In terms of Theorem 1, $G-W-e$ has no $P_{\geq2}$-factor. Hence, $G$ is not $(P_{\geq2},n)$-factor critical avoidable.

\medskip

\noindent{\bf Remark 2.} Next, we explain that the condition on $(n+r+2)$-connected in Theorem 6 is best possible.

Let $G=K_{n+r+1}+((2r+1)K_1\cup K_2)$, where $n$ and $r$ are two nonnegative integers with $n\geq r$. It is obvious that $G$ is $(n+r+1)$-connected and
$I(G)=\frac{n+r+2}{2r+2}=\frac{1}{2}+\frac{n+1}{2(r+1)}>\frac{1}{2}+\frac{n+1}{2(r+2)}=\frac{n+r+3}{2(r+2)}$. Let $W\subseteq V(K_{n+r+1})\subseteq V(G)$
with $|W|=n$ and $e\in E(K_2)$. Then $G-W-e=K_{r+1}+((2r+3)K_1)$. Let $X=V(K_{r+1})\subseteq V(G-W-e)$. Then we derive
$$
i(G-W-e-X)=2r+3>2(r+1)=2|X|.
$$
In terms of Theorem 1, $G-W-e$ has no $P_{\geq2}$-factor. So $G$ is not $(P_{\geq2},n)$-factor critical avoidable.

\medskip

\section{$(P_{\geq3},n)$-factor critical avoidable graphs}

We first verify the following lemma.

\medskip

\noindent{\textbf{Lemma 1.}} Let $n$ and $r$ be two nonnegative integers, let $G$ be an $(n+r+2)$-connected graph, and let $H=G-W-e$ for any
$W\subseteq V(G)$ with $|W|=n$ and any $e\in E(G-W)$. If $sun(H-X)\geq2|X|+1$ for $X\subseteq V(H)$, then $|X|\geq r+2$.

\medskip

\noindent{\it Proof.} If $|X|=0$, then $sun(H)\geq1$.

On the other hand, since $G$ is $(n+r+2)$-connected, $H$ is $(r+1)$-connected. Thus, we have $sun(H)\leq\omega(H)=1$.

Hence, we obtain $sun(H)=1$. Combining this with $H$ being $(r+1)$-connected, $H$ is a sun.

Note that $G$ is $(n+r+2)$-connected, and so $|V(G)|\geq n+r+3$. Thus, $|V(H)|=|V(G)|-n\geq(n+r+3)-n=r+3\geq3$, which implies that $H$ is a big sun.
Hence, there exist at least three vertices with degree 1 in $H$, and so there exists at least one vertex $v$ with $d_{G-W}(v)=1$. Thus, we acquire
$$
d_G(v)\leq d_{G-W}+|W|=n+1\leq n+r+1,
$$
which contradicts that $G$ is $(n+r+2)$-connected. In what follows, we consider $1\leq|X|\leq r+1$.

In light of $sun(H-X)\geq2|X|+1$, we admit
$$
\omega(H-X)\geq sun(H-X)\geq2|X|+1\geq3.
$$
Thus, we derive
$$
\omega(G-W-X)\geq\omega(G-W-X-e)-1=\omega(H-X)-1\geq3-1=2,
$$
Combining this with $|W|=n$ and $1\leq|X|\leq r+1$, we know that $G$ is at most $(n+r+1)$-connected, which contradicts that $G$ is $(n+r+2)$-connected.
Hence, $|X|\geq r+2$. This completes the proof of Lemma 1. \hfill $\Box$

Next, we raise two sufficient conditions using toughness and isolated toughness for graphs being $(P_{\geq3},n)$-factor critical avoidable
graphs, which are the improvements of Theorems 4 and 5.

\medskip

\noindent{\textbf{Theorem 7.}} Let $n$ and $r$ be two nonnegative integers, and let $G$ be an $(n+r+2)$-connected graph. If its toughness
$t(G)>\frac{n+r+2}{2(r+2)}$, then $G$ is $(P_{\geq3},n)$-factor critical avoidable.

\medskip

\noindent{\it Proof.} For a complete graph $G$, Theorem 7 is true. In the following, we assume that $G$ is not a complete graph. Let $H=G-W-e$ for
any $W\subseteq V(G)$ with $|W|=n$ and any $e\in E(G-W)$. It suffices to prove that $H$ admits a $P_{\geq3}$-factor. On the contrary, we assume that
$H$ has no $P_{\geq3}$-factor. In view of Theorem 2, we obtain
$$
sun(H-X)\geq2|X|+1\eqno(1)
$$
for some subset $X$ of $V(H)$.

In view of (1) and Lemma 1, we infer
\begin{eqnarray*}
\omega(G-W-X)&\geq&\omega(G-W-X-e)-1=\omega(H-X)-1\\
&\geq&sun(H-X)-1\geq2|X|+1-1=2|X|\\
&\geq&2(r+2)\geq4.
\end{eqnarray*}
Combining this with the definition of $t(G)$, we get
$$
t(G)\leq\frac{|W\cup X|}{\omega(G-W-X)}\leq\frac{n+|X|}{2|X|}=\frac{1}{2}+\frac{n}{2|X|}\leq\frac{1}{2}+\frac{n}{2(r+2)}=\frac{n+r+2}{2(r+2)},
$$
which contradicts $t(G)>\frac{n+r+2}{2(r+2)}$. We finish the proof of Theorem 7. \hfill $\Box$

\medskip

\noindent{\bf Remark 3.} Next, we show that the condition on $t(G)$ in Theorem 7 is sharp.

Let $G=K_{n+r+2}+((2r+3)K_1\cup K_2)$, where $n$ and $r$ are two nonnegative integers. Obviously, $G$ is $(n+r+2)$-connected and
$t(G)=\frac{n+r+2}{2(r+2)}$. Let $W\subseteq V(K_{n+r+2})\subseteq V(G)$ with $|W|=n$ and $e\in E(K_2)$. Then $G-W-e=K_{r+2}+((2r+5)K_1)$.
Select $X=V(K_{r+2})$ in $G-W-e$. Thus, we derive
$$
sun(G-W-e-X)=2r+5>2(r+2)=2|X|.
$$
In terms of Theorem 2, $G-W-e$ has no $P_{\geq3}$-factor. Hence, $G$ is not $(P_{\geq3},n)$-factor critical avoidable.

\medskip

\noindent{\bf Remark 4.} Next, we explain that the condition on $(n+r+2)$-connected in Theorem 7 cannot be replaced by $(n+r+1)$-connected.

Let $G=K_{n+r+1}+((2r+1)K_1\cup K_2)$, where $n\geq1$ and $r\geq0$ are two integers. We see that $G$ is $(n+r+1)$-connected and
$t(G)=\frac{n+r+1}{2r+2}=\frac{1}{2}+\frac{n}{2(r+1)}>\frac{1}{2}+\frac{n}{2(r+2)}=\frac{n+r+2}{2(r+2)}$. Let $W\subseteq V(K_{n+r+1})\subseteq V(G)$
with $|W|=n$ and $e\in E(K_2)$. Then $G-W-e=K_{r+1}+((2r+3)K_1)$. Choose $X=V(K_{r+1})$ in $G-W-e$. Thus, we admit
$$
sun(G-W-e-X)=2r+3>2(r+1)=2|X|.
$$
In terms of Theorem 2, $G-W-e$ has no $P_{\geq3}$-factor. Therefore, $G$ is not $(P_{\geq3},n)$-factor critical avoidable.

\medskip

\noindent{\textbf{Theorem 8.}} Let $n$ and $r$ be two nonnegative integers, and let $G$ be an $(n+r+2)$-connected graph. If its isolated toughness
$I(G)>\frac{n+3(r+2)}{2(r+2)}$, then $G$ is $(P_{\geq3},n)$-factor critical avoidable.

\medskip

\noindent{\it Proof.} It is obvious that Theorem 8 is true for a complete graph. In what follows, we assume that $G$ is not complete. Let
$H=G-W-e$ for any $W\subseteq V(G)$ with $|W|=n$ and any $e=uv\in E(G-W)$. It suffices to verify that $H$ contains a $P_{\geq3}$-factor. By
means of contrary, we assume that $H$ has no $P_{\geq3}$-factor. Then by Theorem 2, we admit
$$
sun(H-X)\geq2|X|+1\eqno(1)
$$
for some vertex subset $X$ of $H$.

Suppose that there exist $a$ isolated vertices, $b$ $K_2$'s and $c$ big sun components $H_1,H_2,\cdots,H_c$, where $|V(H_i)|\geq6$, in $H-X$.
Let $R_i$ be the factor-critical subgraph of $H_i$. We select one vertex from every $K_2$ component of $H-X$, and denote the set of such
vertices by $Y$. Thus, we admit
$$
sun(H-X)=a+b+c.\eqno(2)
$$

According to (1), (2), Lemma 1 and $|V(R_i)|\geq3$, we infer
$$
a+b+\sum\limits_{i=1}^{c}{|V(R_i)|}\geq a+b+3c\geq a+b+c=sun(H-X)\geq2|X|+1\geq2(r+2)+1\geq5.\eqno(3)
$$

In fact, $sun(G-W-X-e)+1\geq sun(G-W-X)\geq sun(G-W-X-e)-2$. The following proof is divided into four cases.

\noindent{\bf Case 1.} $sun(G-W-X)=sun(G-W-X-e)-2$.

In this case, $u\in V(aK_1)$ and $v\in V(bK_2)$, or $u\in V(aK_1)$ and $v\in V(H_i)$, or $u,v$ belong to two different $K_2$ components, or
$u\in V(H_i)$ and $v\in V(bK_2)$, or $v\in V(H_i)$ and $u\in V(H_j)$ ($i\neq j$).

\noindent{\bf Claim 1.} $I(G)\leq\frac{n+|X|+b+\sum\limits_{i=1}^{c}{|V(R_i)|}}{a+b+\sum\limits_{i=1}^{c}{|V(R_i)|}}$.

\noindent{\it Proof.} We consider the following two subcases.

\noindent{\bf Subcase 1.1.} $u\in V(aK_1)$ and $v\in V(bK_2)$, or $u\in V(H_i)$ and $v\in V(bK_2)$, or $u,v$ belong to two different $K_2$ components.

We choose such $Y$ with $v\in Y$. Then we deduce
\begin{eqnarray*}
i(G-W-X-Y-\bigcup\limits_{i=1}^{c}{V(R_i)})&=&i(G-W-X-e-Y-\bigcup\limits_{i=1}^{c}{V(R_i)})\\
&=&i(H-X-Y-\bigcup\limits_{i=1}^{c}{V(R_i)})\\
&=&a+b+\sum\limits_{i=1}^{c}{|V(R_i)|}.
\end{eqnarray*}
In terms of (3) and the definition of $I(G)$, we get
$$
I(G)\leq\frac{|W\cup X\cup Y\cup(\bigcup\limits_{i=1}^{c}{V(R_i)})|}{i(G-W-X-Y-\bigcup\limits_{i=1}^{c}{V(R_i)})}
=\frac{n+|X|+b+\sum\limits_{i=1}^{c}{|V(R_i)|}}{a+b+\sum\limits_{i=1}^{c}{|V(R_i)|}}.
$$

\noindent{\bf Subcase 1.2.} $u\in V(aK_1)$ and $v\in V(H_i)$, or $v\in V(H_i)$ and $u\in V(H_j)$ ($i\neq j$).

If $v\in V(R_i)$, then we write $Z=\bigcup\limits_{i=1}^{c}{V(R_i)}$. If $v\in V(H_i)\setminus V(R_i)$, then $d_{H_i}(v)=1$, which implies
that there exists $w\in V(R_i)$ such that $vw\in E(H_i)$. Now we write $Z=((\bigcup\limits_{i=1}^{c}{V(R_i)})\cup\{v\})\setminus\{w\}$.
Thus, we refer
\begin{eqnarray*}
i(G-W-X-Y-Z)&=&i(G-W-X-e-Y-Z)\\
&=&i(H-X-Y-Z)\\
&=&a+b+\sum\limits_{i=1}^{c}{|V(R_i)|}.
\end{eqnarray*}
Using (3) and the definition of $I(G)$, we derive
$$
I(G)\leq\frac{|W\cup X\cup Y\cup Z|}{i(G-W-X-Y-Z)}
=\frac{n+|X|+b+\sum\limits_{i=1}^{c}{|V(R_i)|}}{a+b+\sum\limits_{i=1}^{c}{|V(R_i)|}}.
$$
This completes the proof of Claim 1. \hfill $\Box$

In light of (3), Claim 1 and Lemma 1, we obtain
\begin{eqnarray*}
I(G)&\leq&\frac{n+|X|+b+\sum\limits_{i=1}^{c}{|V(R_i)|}}{a+b+\sum\limits_{i=1}^{c}{|V(R_i)|}}\\
&\leq&\frac{n+|X|+a+b+\sum\limits_{i=1}^{c}{|V(R_i)|}}{a+b+\sum\limits_{i=1}^{c}{|V(R_i)|}}\\
&=&1+\frac{n+|X|}{a+b+\sum\limits_{i=1}^{c}{|V(R_i)|}}\\
&\leq&1+\frac{n+|X|}{2|X|+1}\\
&<&1+\frac{n+|X|+\frac{1}{2}}{2|X|+1}\\
&=&\frac{3}{2}+\frac{n}{2|X|+1}\\
&<&\frac{3}{2}+\frac{n}{2|X|}\\
&\leq&\frac{3}{2}+\frac{n}{2(r+2)}\\
&=&\frac{n+3(r+2)}{2(r+2)},
\end{eqnarray*}
which contradicts $I(G)>\frac{n+3(r+2)}{2(r+2)}$.

\noindent{\bf Case 2.} $sun(G-W-X)=sun(G-W-X-e)-1$.

In this case, $u,v\in V(aK_1)$, or $u\in V(M)$ and $v\in V(Q)$, where $M$ is a sun component of $H-X$, $Q$ is a non-sun component of $H-X$,
and $M\cup Q\cup\{e\}$ is a non-sun component of $(H-X)\cup\{e\}$.

\noindent{\bf Subcase 2.1.} $u,v\in V(aK_1)$.

In this subcase, $a\geq2$. Thus, we admit
\begin{eqnarray*}
i(G-W-X-Y-\bigcup\limits_{i=1}^{c}{V(R_i)}-v)&=&i(G-W-X-e-Y-\bigcup\limits_{i=1}^{c}{V(R_i)}-v)\\
&=&i(H-X-Y-\bigcup\limits_{i=1}^{c}{V(R_i)}-v)\\
&=&a+b+\sum\limits_{i=1}^{c}{|V(R_i)|}-1.
\end{eqnarray*}
Combining this with (3), $a\geq2$, Lemma 1 and the definition of $I(G)$, we derive
\begin{eqnarray*}
I(G)&\leq&\frac{|W\cup X\cup Y\cup(\bigcup\limits_{i=1}^{c}{V(R_i)})\cup\{v\}|}{i(G-W-X-Y-\bigcup\limits_{i=1}^{c}{V(R_i)}-v)}\\
&=&\frac{n+|X|+b+\sum\limits_{i=1}^{c}{|V(R_i)|}+1}{a+b+\sum\limits_{i=1}^{c}{|V(R_i)|}-1}\\
&=&1+\frac{n+|X|+2-a}{a+b+\sum\limits_{i=1}^{c}{|V(R_i)|}-1}\\
&\leq&1+\frac{n+|X|}{2|X|}\\
&=&\frac{3}{2}+\frac{n}{2|X|}\\
&\leq&\frac{3}{2}+\frac{n}{2(r+2)}\\
&=&\frac{n+3(r+2)}{2(r+2)},
\end{eqnarray*}
which contradicts $I(G)>\frac{n+3(r+2)}{2(r+2)}$.

\noindent{\bf Subcase 2.2.} $u\in V(M)$ and $v\in V(Q)$, where $M$ is a sun component of $H-X$, $Q$ is a non-sun component of $H-X$,
and $M\cup Q\cup\{e\}$ is a non-sun component of $(H-X)\cup\{e\}$.

If $M=K_1$, then $a\geq1$. Thus, we possess
\begin{eqnarray*}
i(G-W-X-Y-\bigcup\limits_{i=1}^{c}{V(R_i)}-v)&=&i(G-W-X-e-Y-\bigcup\limits_{i=1}^{c}{V(R_i)}-v)\\
&=&i(H-X-Y-\bigcup\limits_{i=1}^{c}{V(R_i)}-v)\\
&=&a+b+\sum\limits_{i=1}^{c}{|V(R_i)|}.
\end{eqnarray*}
Then using (3), $a\geq1$, Lemma 1, the definition of $I(G)$ and $I(G)>\frac{n+3(r+2)}{2(r+2)}$, we infer
\begin{eqnarray*}
\frac{n+3(r+2)}{2(r+2)}&<&I(G)\leq\frac{|W\cup X\cup Y\cup(\bigcup\limits_{i=1}^{c}{V(R_i)})\cup\{v\}|}{i(G-W-X-Y-\bigcup\limits_{i=1}^{c}{V(R_i)}-v)}\\
&=&\frac{n+|X|+b+\sum\limits_{i=1}^{c}{|V(R_i)|}+1}{a+b+\sum\limits_{i=1}^{c}{|V(R_i)|}}\\
&=&1+\frac{n+|X|+1-a}{a+b+\sum\limits_{i=1}^{c}{|V(R_i)|}}\\
&\leq&1+\frac{n+|X|}{2|X|+1}\\
&<&1+\frac{n+|X|}{2|X|}\\
&=&\frac{3}{2}+\frac{n}{2|X|}\\
&\leq&\frac{3}{2}+\frac{n}{2(r+2)}\\
&=&\frac{n+3(r+2)}{2(r+2)},
\end{eqnarray*}
which is a contradiction.

If $M=K_2$, then we choose such $Y$ with $u\in Y$ and select $Z=\bigcup\limits_{i=1}^{c}{V(R_i)}$. If $M=H_i$, then $u\in V(R_i)$ or
$u\in V(H_i)\setminus V(R_i)$. If $u\in V(R_i)$, then we choose $Z=\bigcup\limits_{i=1}^{c}{V(R_i)}$. If $u\in V(H_i)\setminus V(R_i)$, then there
exists $w\in V(R_i)$ such that $uw\in E(H_i)$. Thus, we select $Z=((\bigcup\limits_{i=1}^{c}{V(R_i)})\cup\{u\})\setminus\{w\}$.

Hence, we deduce
\begin{eqnarray*}
i(G-W-X-Y-Z)&=&i(G-W-X-e-Y-Z)\\
&=&i(H-X-Y-Z)\\
&=&a+b+\sum\limits_{i=1}^{c}{|V(R_i)|}.
\end{eqnarray*}
Combining this with (3), Lemma 1 and then definition of $I(G)$, we get
\begin{eqnarray*}
I(G)&\leq&\frac{|W\cup X\cup Y\cup Z|}{i(G-W-X-Y-Z)}\\
&=&\frac{n+|X|+b+\sum\limits_{i=1}^{c}{|V(R_i)|}}{a+b+\sum\limits_{i=1}^{c}{|V(R_i)|}}\\
&\leq&\frac{n+|X|+a+b+\sum\limits_{i=1}^{c}{|V(R_i)|}}{a+b+\sum\limits_{i=1}^{c}{|V(R_i)|}}\\
&=&1+\frac{n+|X|}{a+b+\sum\limits_{i=1}^{c}{|V(R_i)|}}\\
&\leq&1+\frac{n+|X|}{2|X|+1}\\
&<&1+\frac{n+|X|}{2|X|}\\
&=&\frac{3}{2}+\frac{n}{2|X|}\\
&\leq&\frac{3}{2}+\frac{n}{2(r+2)}\\
&=&\frac{n+3(r+2)}{2(r+2)},
\end{eqnarray*}
which contradicts $I(G)>\frac{n+3(r+2)}{2(r+2)}$.

\noindent{\bf Case 3.} $sun(G-W-X)=sun(G-W-X-e)$.

In this case, $u\in V(R_i)$ and $v\in V(H_i)\setminus V(R_i)$, or $u,v\in V(M)$, or $u\in V(M)$ and $v\in V(Q)$, where $M$ is a non-sun component of
$H-X$, $M\cup\{e\}$ is a non-sun component of $(H-X)\cup\{e\}$, $Q$ is a non-sun component of $H-X$, and $M\cup Q\cup\{e\}$ is a non-sun component
of $(H-X)\cup\{e\}$. Then we infer
\begin{eqnarray*}
i(G-W-X-Y-\bigcup\limits_{i=1}^{c}{V(R_i)})&=&i(G-W-X-e-Y-\bigcup\limits_{i=1}^{c}{V(R_i)})\\
&=&i(H-X-Y-\bigcup\limits_{i=1}^{c}{V(R_i)})\\
&=&a+b+\sum\limits_{i=1}^{c}{|V(R_i)|}.
\end{eqnarray*}
It follows from (3), Lemma 1, the definition of $I(G)$ and $I(G)>\frac{n+3(r+2)}{2(r+2)}$ that
\begin{eqnarray*}
\frac{n+3(r+2)}{2(r+2)}&<I(G)&\leq\frac{|W\cup X\cup Y\cup(\bigcup\limits_{i=1}^{c}{V(R_i)})|}{i(G-W-X-Y-\bigcup\limits_{i=1}^{c}{V(R_i)})}\\
&=&\frac{n+|X|+b+\sum\limits_{i=1}^{c}{|V(R_i)|}}{a+b+\sum\limits_{i=1}^{c}{|V(R_i)|}}\\
&\leq&1+\frac{n+|X|}{a+b+\sum\limits_{i=1}^{c}{|V(R_i)|}}\\
&\leq&1+\frac{n+|X|}{2|X|+1}\\
&<&1+\frac{n+|X|}{2|X|}\\
&=&\frac{3}{2}+\frac{n}{2|X|}\\
&\leq&\frac{3}{2}+\frac{n}{2(r+2)}\\
&=&\frac{n+3(r+2)}{2(r+2)},
\end{eqnarray*}
which is a contradiction.

\noindent{\bf Case 4.} $sun(G-W-X)=sun(G-W-X-e)+1$.

In this case, $u,v\in V(R_i)$. The following proof is similar to that of Case 3, we easily deduce
$$
I(G)\leq\frac{|W\cup X\cup Y\cup(\bigcup\limits_{i=1}^{c}{V(R_i)})|}{i(G-W-X-Y-\bigcup\limits_{i=1}^{c}{V(R_i)})}<\frac{n+3(r+2)}{2(r+2)},
$$
which contradicts $I(G)>\frac{n+3(r+2)}{2(r+2)}$. We finish the proof of Theorem 8. \hfill $\Box$

\medskip

\noindent{\bf Remark 5.} Now, we show that the condition $I(G)>\frac{n+3(r+2)}{2(r+2)}$ in Theorem 8 cannot be replaced by
$I(G)\geq\frac{n+3(r+2)}{2(r+2)}$.

Let $G=K_{n+r+2}+((2r+4)K_2)$, where $n$ and $r$ are two nonnegative integers. It is obvious that $G$ is $(n+r+2)$-connected and
$I(G)=\frac{n+3(r+2)}{2(r+2)}$. Set $W\subseteq V(K_{n+r+2})\subseteq V(G)$ with $|W|=n$ and $e\in E(K_2)$. Then $G-W-e=K_{r+2}+((2r+3)K_2\cup (2K_1))$.
Let $X=V(K_{r+2})$ in $G-W-e$. Thus, we acquire
$$
sun(G-W-e-X)=2r+5>2(r+2)=2|X|.
$$
According to Theorem 2, $G-W-e$ has no $P_{\geq3}$-factor. So $G$ is not $(P_{\geq3},n)$-factor critical avoidable.

\medskip

\noindent{\bf Remark 6.} Now, we claim that the condition that $(n+r+2)$-connected in Theorem 8 is sharp.

Let $G=K_{n+r+1}+((2r+2)K_2)$, where $n\geq1$ and $r\geq0$ are two integers. Clearly, $G$ is $(n+r+1)$-connected and
$I(G)=\frac{n+3(r+1)}{2(r+1)}=\frac{3}{2}+\frac{n}{2(r+1)}>\frac{3}{2}+\frac{n}{2(r+2)}=\frac{n+3(r+2)}{2(r+2)}$. Let $W\subseteq V(K_{n+r+1})\subseteq V(G)$
with $|W|=n$ and $e\in E(K_2)$. Then $G-W-e=K_{r+1}+((2r+1)K_1\cup(2K_1))$. Let $X=V(K_{r+1})$ in $G-W-e$. Thus, we refer
$$
sun(G-W-e-X)=2r+3>2(r+1)=2|X|.
$$
Using Theorem 2, $G-W-e$ has no $P_{\geq3}$-factor. Hence, $G$ is not $(P_{\geq3},n)$-factor critical avoidable.

\medskip

%\section*{Acknowledgments}
%This work was supported by the National Natural Science Foundation of China (Nos. 11731002).

\end{document}